\documentclass[12pt]{article}
\usepackage{amssymb,amsmath}
\usepackage{cases}
\usepackage{amsfonts}
\usepackage{color,xcolor}
\usepackage[left=2.0cm,right=2.0cm,top=2.0cm,bottom=2.0cm]{geometry}
\usepackage[colorlinks,citecolor=blue,urlcolor=blue]{hyperref}
\usepackage[utf8]{inputenc}

\newtheorem{theorem}{Theorem}[section]

\newtheorem{lemma}{Lemma}[section]
\newtheorem{proposition}{Proposition}[section]

\newtheorem{definition}{Definition}[section]
\newtheorem{remark}{Remark}[section]

\newcommand{\bal}{\begin{align}}
\newcommand{\bbal}{\begin{align*}}
\newcommand{\beq}{\begin{equation}}
\newcommand{\eeq}{\end{equation}}
\newcommand{\bca}{\begin{cases}}
\newcommand{\eca}{\end{cases}}

\newcommand{\pa}{\partial}
\newcommand{\fr}{\frac}

\newcommand{\dd}{\mathrm{d}}

\newcommand{\R}{\mathbb{R}}

\newcommand{\les}{\lesssim}

\newcommand{\bi}{\Big}

\begin{document}
\title{Well-posedness and Continuity Properties of the Degasperis-Procesi equation in critical Besov space}

\author{Jinlu Li$^{1}$, Yanghai Yu$^{2,}$\footnote{E-mail: lijinlu@gnnu.edu.cn; yuyanghai214@sina.com(Corresponding author); mathzwp2010@163.com} and Weipeng Zhu$^{3}$\\
\small $^1$ School of Mathematics and Computer Sciences, Gannan Normal University, Ganzhou 341000, China\\
\small $^2$ School of Mathematics and Statistics, Anhui Normal University, Wuhu 241002, China\\
\small $^3$ School of Mathematics and Big Data, Foshan University, Foshan, Guangdong 528000, China}

\date{\today}

\maketitle\noindent{\hrulefill}

{\bf Abstract:} In this paper, we obtain the local-in-time existence and uniqueness of solution to the Degasperis-Procesi equation in $B^1_{\infty,1}(\R)$. Moreover, we prove that the data-to-solution of this equation is continuous but not uniformly continuous in $B^1_{\infty,1}$.


\vskip0mm\noindent{\hrulefill}

\section{Introduction}\label{sec1}
In this paper, we consider the Cauchy problem for the Degasperis-Procesi (DP) equation \cite{DP1}
\begin{align}\label{dp}
\begin{cases}
\pa_tu+uu_x=-\frac32\pa_x(1-\pa^2_x)^{-1}(u^2), &\quad (t,x)\in \R^+\times\R,\\
u(0,x)=u_0(x), &\quad x\in \R.
\end{cases}
\end{align}
We should also mention the following famous Camassa-Holm (CH) equation \cite{Camassa}
\begin{equation}\label{CH}
\begin{cases}
\partial_tu+u\pa_xu=-\pa_x\big(1-\pa^2_x\big)^{-1}\big(u^2+\fr12(\pa_xu)^2\big), \; &(x,t)\in \R\times\R^+,\\
u(0,x)=u_0(x),\; &x\in \R.
\end{cases}
\end{equation}
The DP equation can be regarded as a model for nonlinear shallow water dynamics and its asymptotic
accuracy is the same as for the CH shallow water equation \cite{Dullin}. The DP equation is integrable and has a
bi-Hamiltonian structure \cite{DP}. An inverse scattering approach for computing $n$-peakon solutions to the DP
equation was presented in \cite{Lundmark}. Its traveling wave solutions were investigated in \cite{Lenells,Vakhnenko}.
Regarding well-posedness (existence, uniqueness, and stability of solutions) of the Cauchy problem for the
DP equation \eqref{dp}, we refer to see \cite{Gui,Himonas,Li0,Lin,Yin2}. Similar to the CH equation, the DP equation has also global strong solutions
\cite{Liu2006,Yin2003,Yin2004} and finite time blow-up solutions \cite{Escherjfa,Escher,Liu2006,Liuyin,Yin2}.
Although the DP equation is similar to the CH equation in several aspects, these two equations are truly
different. One of the novel features of the DP different from the CH equation is that it has not only peakon solutions \cite{DP} and periodic peakon solutions \cite{YinJFA}, but also shock peakons \cite{Lundmark2007} and the periodic shock waves \cite{Escher}.

In our recent papers \cite{Li1,Li2}, we proved that the DP equation \eqref{dp} is not uniformly continuous in $B^s_{p,r}$ with $s>1+1/p, \ p\in [1,\infty], \ r\in [1,\infty)$ or $s=1+1/p, \ p\in [1,\infty),\  r=1.$ Guo, Liu, Molinet and Yin \cite{Guo-Yin} established the ill-posedness of \eqref{dp} in the critical Sobolev space $H^{3/2}$ and even in the Besov space $B_{p,r}^{1+1/p}$ with $p\in[1,\infty],r\in(1,\infty]$ by proving the norm inflation. For more ill-posedness results of the Cauchy problem for the DP equation \eqref{dp}, we refer to see \cite{Himonas2014,Himonas2016} and the references therein. However, to our best knowledge, for the Cauchy problem of the DP equation \eqref{dp} with
initial data in $B^{1}_{\infty,1}$ is locally well-posed in the sense of Hadamard has not been solved yet. Before stating our main result precisely, we recall the notion of well-posedness in the sense of Hadamard:
\begin{definition}[Local Well-posedness] We say that the Cauchy problem \eqref{dp} is locally well-posed
in a Banach space $X$ if the following three conditions hold
\begin{enumerate}
  \item (Local existence)\; For any initial data $u_0\in X$, there exists a short time $T = T(u_0) > 0$ and a solution $\mathbf{S}_{t}(u_0)\in\mathcal{C}([0,T),X)$ to the Cauchy problem \eqref{dp};
  \item (Uniqueness)\; This solution $\mathbf{S}_{t}(u_0)$ is unique in the space $\mathcal{C}([0,T),X)$;
  \item (Continuous Dependence)\; The data-to-solution map $u_0 \mapsto \mathbf{S}_{t}(u_0)$ is continuous in the following
sense:  for any $T_1 < T$ and $\varepsilon > 0$, there exists $ \delta> 0$, such that if $\|u_0-\widetilde{u}_0\|_{X}\leq \delta$,
then $\mathbf{S}_{t}(\tilde{u}_0)$ exists up to $T_1$ and
$$\|\mathbf{S}_{t}(u_0)-\mathbf{S}_{t}(\tilde{u}_0)\|_{\mathcal{C}([0,T),X)}\leq \varepsilon.$$
\end{enumerate}
\end{definition}
Now, we can formulate the main result.
\begin{theorem}\label{th1}
Let $u_0\in B^{1}_{\infty,1}(\R)$. Then, there exists some time $T>0$ such that
\begin{enumerate}
  \item system \eqref{dp} has a solution $u\in \mathcal{C}([0,T],B^{1}_{\infty,1})\cap \mathcal{C}^1([0,T],B^{0}_{\infty,1})$;
  \item the solutions of \eqref{dp} are unique;
  \item if the initial data have some additional regularity $u_0\in B^{s}_{\infty,1}$ with $s>1$, then the solution exists on $[0, T]$ and is in
the space $\mathcal{C}([0,T],B^{s}_{\infty,1})\cap \mathcal{C}^1([0,T],B^{s-1}_{\infty,1})$;
  \item the data-to-solution map $u_0 \mapsto u(t)$ is continuous from any bounded
subset of $u_0\in B^{1}_{\infty,1}$ into $\mathcal{C}([0,T],B^{1}_{\infty,1})$;
\item the data-to-solution map $u_0 \mapsto u(t)$ is not uniformly continuous from any bounded subset of $u_0\in B^{1}_{\infty,1}$ into $\mathcal{C}([0,T],B^{1}_{\infty,1})$.
\end{enumerate}
\end{theorem}
\begin{remark}\label{re0}
Theorem \ref{th1} tells us that there exists a positive time $T$ such that System \eqref{dp} has a unique solution and the solution is continuously dependent on the initial data in the sense of Hadmard.
\end{remark}
\section{Preliminaries}\label{sec2}
{\bf Notation}\quad In the following, we denote by $*$ the convolution. Throughout this paper, we shall use
$\big(1-\pa^2_x\big)^{-1}f=p*f$ with $p(x)=\fr12e^{-|x|}.$
Given a Banach space $X$, we denote its norm by $\|\cdot\|_{X}$. For $I\subset\R$, we denote by $\mathcal{C}(I;X)$ the set of continuous functions on $I$ with values in $X$. Sometimes we will denote $L^p(0,T;X)$ by $L_T^pX$.
For all $f\in \mathcal{S}'$, the Fourier transform $\mathcal{F}f$ (also denoted by $\widehat{f}$) is defined by
$$
\mathcal{F}f(\xi)=\widehat{f}(\xi)=\int_{\R}e^{-ix\xi}f(x)\dd x \quad\text{for any}\; \xi\in\R.
$$
Next, we will recall some facts about the Littlewood-Paley decomposition, the nonhomogeneous Besov spaces and their some useful properties.
\begin{definition}[See \cite{B.C.D}] Let $\mathcal{B}:=\{\xi\in\mathbb{R}:|\xi|\leq 4/3\}$ and $\mathcal{C}:=\{\xi\in\mathbb{R}:3/4\leq|\xi|\leq 8/3\}.$
There exist two radial functions $\chi\in C_c^{\infty}(\mathcal{B})$ and $\varphi\in C_c^{\infty}(\mathcal{C})$ both taking values in $[0,1]$ such that
$$
\chi(\xi)+\sum_{j\geq0}\varphi(2^{-j}\xi)=1 \quad \forall \;  \xi\in \R.$$
For every $u\in \mathcal{S'}(\mathbb{R})$, the Littlewood-Paley dyadic blocks ${\Delta}_j$ are defined as follows
\begin{numcases}{\Delta_ju=}
0, & if $j\leq-2$;\nonumber\\
\chi(D)u=\mathcal{F}^{-1}(\chi \mathcal{F}u), & if $j=-1$;\nonumber\\
\varphi(2^{-j}D)u=\mathcal{F}^{-1}\big(\varphi(2^{-j}\cdot)\mathcal{F}u\big), & if $j\geq0$.\nonumber
\end{numcases}
The inhomogeneous low-frequency cut-off operator $S_j$ is defined by
$$
S_ju=\sum_{q=-1}^{j-1}{\Delta}_qu.
$$

\end{definition}
\begin{definition}[See \cite{B.C.D}]
Let $s\in\mathbb{R}$ and $(p,r)\in[1, \infty]^2$. The nonhomogeneous Besov space $B^{s}_{p,r}(\R)$ is defined by
\begin{align*}
B^{s}_{p,r}(\R):=\Big\{f\in \mathcal{S}'(\R):\;\|f\|_{B^{s}_{p,r}(\mathbb{R})}<\infty\Big\},
\end{align*}
where
\begin{numcases}{\|f\|_{B^{s}_{p,r}(\mathbb{R})}=}
\left(\sum_{j\geq-1}2^{sjr}\|\Delta_jf\|^r_{L^p(\mathbb{R})}\right)^{\fr1r}, &if $1\leq r<\infty$,\nonumber\\
\sup_{j\geq-1}2^{sj}\|\Delta_jf\|_{L^p(\mathbb{R})}, &if $r=\infty$.\nonumber
\end{numcases}
\end{definition}
\begin{lemma}[See \cite{B.C.D}]\label{le2}
For $p\in[1, \infty]$ and $s>0$, $B^s_{p,1}(\R)\cap L^\infty(\R)$ is an algebra. Moreover, $B^{0}_{\infty,1}(\R)\hookrightarrow L^\infty(\R)\hookrightarrow B^{0}_{\infty,\infty}(\R)$, and for any $u,v \in B^s_{p,1}(\R)\cap L^\infty(\R)$, we have
\bbal
&\|uv\|_{B^{s}_{p,1}}\leq C(\|u\|_{B^{s}_{p,1}}\|v\|_{L^\infty}+\|v\|_{B^{s}_{p,1}}\|u\|_{L^\infty}).
\end{align*}
We also have the following interpolation inequality
\bbal
&\|u\|_{B^{1}_{\infty,1}}\leq C\|u\|_{B^{0}_{\infty,\infty}}^{\fr12}\|u\|_{B^{2}_{\infty,\infty}}^{\fr12}.
\end{align*}
\end{lemma}
\begin{lemma}[See \cite{B.C.D,Li0}]\label{le3}
Let $(p,r)\in[1, \infty]^2$ and $\sigma\geq-\min\big\{\frac1p, 1-\frac1p\big\}$. Assume that $f_0\in B^\sigma_{p,r}(\R)$, $g\in L^1([0,T]; B^\sigma_{p,r}(\R))$ and $\pa_{x}\mathbf{u}\in
L^1([0,T]; B^{\sigma-1}_{p,r}(\R))$ if $\sigma>1+\fr1p$  or $\sigma=1+\fr1p, r=1$.
If $f\in L^\infty([0,T]; B^\sigma_{p,r}(\R))\cap \mathcal{C}([0,T]; \mathcal{S}'(\R))$ solves the following linear transport equation:
\begin{equation*}
\quad \partial_t f+\mathbf{u}\pa_xf=g,\quad \; f|_{t=0} =f_0.
\end{equation*}
 Then there exists a constant $C=C(p,r,\sigma)$ such that the following statement holds
\begin{align*}
\|f(t)\|_{B^\sigma_{p,r}}\leq e^{CV(t)} \Big(\|f_0\|_{B^\sigma_{p,r}}+\int_0^t e^{-CV(\tau)} \|g(\tau)\|_{B^\sigma_{p,r}}\mathrm{d}\tau\Big),
\end{align*}
where
$$V(t)=\int_0^t \|\pa_x\mathbf{u}(\tau)\|_{B^{\sigma-1}_{p,r}}\mathrm{d}\tau\quad \text{if}\quad \sigma>1+{\fr 1p} \quad \text{or} \quad \{\sigma=1+\fr1p,\; r=1\}.$$
\end{lemma}

\section{Proof of Theorem \ref{th1}}
In this section, we divide Proof of Theorem \ref{th1} into three steps.
\subsection{Existence and Uniqueness}
Applying Lemma \ref{le3} to Eq.\eqref{dp}, we have for $t\in[0,T]$
\bal\label{jl}
\|u(t)\|_{B^{1}_{\infty,1}}&\leq e^{V(t)}\Big(\|u_0\|_{B^{1}_{\infty,1}}+\int_0^te^{-V(\tau)}\|\mathbf{P}(u)(\tau)\|_{B^{1}_{\infty,1}}\dd \tau\Big)\nonumber\\
&\leq e^{V(t)}\Big(\|u_0\|_{B^{1}_{\infty,1}}+\int_0^te^{-V(\tau)}\|u(\tau)\|^2_{B^{1}_{\infty,1}}\dd \tau\Big),
\end{align}
where we denote
$$V(t)=\int_0^t \|u\|_{B^{1}_{\infty,1}}\dd\tau$$
and have used
$$\|\mathbf{P}(u)\|_{B^{1}_{\infty,1}}\leq C\|u\|^2_{B^{1}_{\infty,1}}.$$
Then we obtain form \eqref{jl} that
\bbal
F(t):=e^{-V(t)}\|u(t)\|_{B^{1}_{\infty,1}}&\leq \|u_0\|_{B^{1}_{\infty,1}}+\int_0^tF(\tau)\|u(\tau)\|_{B^{1}_{\infty,1}}\dd \tau,
\end{align*}
which follows from the Gronwall inequality that
\bbal
\|u(t)\|_{B^{1}_{\infty,1}}&\leq \|u_0\|_{B^{1}_{\infty,1}}\exp\Big(C\int_0^t \|u\|_{B^{1}_{\infty,1}}\dd\tau\Big).
\end{align*}
By the continuity argument, we can obtain the uniform bounds
\bal\label{l}
\|u\|_{L^\infty_TB^{1}_{\infty,1}}&\leq C\|u_0\|_{B^{1}_{\infty,1}}.
\end{align}
Applying Lemma \ref{le3} to Eq.\eqref{dp} once again, we have for all $t\in[0,T]$ and $s>1$
\bal\label{u2}
\|u(t)\|_{L^\infty_TB^{s}_{\infty,1}}\leq \|u_0\|_{B^{s}_{\infty,1}}\exp\Big(C\int_0^T \|u\|_{B^{1}_{\infty,1}}\dd\tau\Big)\leq C\|u_0\|_{B^s_{\infty,1}}.
\end{align}
The existence follows the standard procedure, we omit the details. The uniqueness is the direct result of the following lemma. In fact, suppose that $u_1,u_2\in \mathcal{C}([0,T],B^{1}_{\infty,1})$ are two solutions of \eqref{dp} with the same initial data $u_0$, then
we have
$$\|u_1-u_2\|_{L^\infty}\leq C\|u_1(0)-u_2(0)\|_{L^\infty}=0,$$
which implies the uniqueness.
\begin{lemma}\label{ley1} Let $u,v\in \mathcal{C}([0,T],B^{1}_{\infty,1})$ be two solutions of \eqref{dp} associated with $u_0$ and $v_0$, respectively. Then we have the estimate for the difference $w=u-v$
\begin{align}
&\|w\|_{L^\infty}\leq C\|w_0\|_{L^\infty},\label{y1}
\\&\|w\|_{B^1_{\infty,1}}\leq C\Big(\|w_0\|_{B^1_{\infty,1}}+\int^t_0\|\pa_{x}v\|_{B^1_{\infty,1}}\|w\|_{L^\infty}\dd \tau\Big),\label{y3}
\end{align}
where the constants $C$ depends on $T$ and initial norm $\|u_0,v_0\|_{B^1_{\infty,1}}$.
\end{lemma}
{\bf Proof.}\quad
Obviously, we have
\begin{align}\label{m}
\begin{cases}
\pa_tw+uw_x=-wv_x-\frac32\pa_x(1-\pa^2_x)^{-1}\big(w(u+v)\big), \\
w(0,x)=u_0(x)-v_0(x).
\end{cases}
\end{align}
Taking the inner product of Eq. \eqref{m} with $|w|^{p-2}w$ with $p\geq2$, we obtain
\begin{align}\label{z}
\fr1p\frac{\dd}{\dd t}\|w\|^p_{L^p}&=-\int_{\R}v_x|w|^{p}\dd x+\fr1p\int_{\R}u_x|w|^{p}\dd x-\frac32\int_{\R}\pa_x(1-\pa^2_x)^{-1}\big(w(u+v)\big)|w|^{p-2}w\dd x\nonumber\\
&\leq \|u_x,v_x\|_{L^\infty}\|w\|^p_{L^p}+\frac32\big\|\pa_x(1-\pa^2_x)^{-1}\big(w(u+v)\big)\big\|_{L^p}\|w\|^{p-1}_{L^p}\nonumber\\
&\leq \|u_x,v_x\|_{L^\infty}\|w\|^p_{L^p}+\frac32\|w(u+v)\|_{L^p}\|w\|^{p-1}_{L^p}\nonumber\\
&\leq C\|u,v,u_x,v_x\|_{L^\infty}\|w\|^p_{L^p},
\end{align}
where we have used that
$$\big\|\pa_x(1-\pa^2_x)^{-1}\big(w(u+v)\big)\big\|_{L^p}=\big\|\pa_xp\ast\big(w(u+v)\big)\big\|_{L^p}\leq \|w(u+v)\|_{L^p}.$$
Integrating the above \eqref{z} with respect to time $t\in[0,T]$ yields
\begin{align*}
\|w\|_{L^p}
&\leq \|w_0\|_{L^p}\exp\Big(C\int_0^T\|u,v\|_{B^1_{\infty,1}}\dd s\Big).
\end{align*}
Letting $p=\infty$ implies \eqref{y1}.

Applying Lemma \ref{le3} to Eq \eqref{m} yields
\begin{align*}
\|w(t)\|_{B^1_{\infty,1}}-C\|w_0\|_{B^1_{\infty,1}}&\les\int_0^t \|wv_x,\frac32\pa_x(1-\pa^2_x)^{-1}\big(w(u+v)\big)\|_{B^1_{\infty,1}}\mathrm{d}\tau\\
&\les \int_0^t\|wv_x\|_{B^1_{\infty,1}}+\|\pa_x(1-\pa^2_x)^{-1}\big(w(u+v)\big)\|_{B^1_{\infty,1}}\mathrm{d}\tau\\
&\les\int_0^t\|wv_x\|_{B^1_{\infty,1}}+\|w(u+v)\|_{B^1_{\infty,1}}\mathrm{d}\tau\\
&\les\int_0^t\|w\|_{B^1_{\infty,1}}\|v_x\|_{L^\infty}+\|v_x\|_{B^1_{\infty,1}}\|w\|_{L^\infty}+\|w\|_{B^1_{\infty,1}}\|u+v\|_{B^1_{\infty,1}}\mathrm{d}\tau\\
&\les\int_0^t\|w\|_{B^1_{\infty,1}}\|u,v\|_{B^1_{\infty,1}}\mathrm{d}\tau+\int_0^t\|v_x\|_{B^1_{\infty,1}}\|w\|_{L^\infty}\mathrm{d}\tau.
\end{align*}
Gronwall's inequality gives us the desired \eqref{y3}.

\subsection{Continuous Dependence}
Now we will prove that the solution is continuously dependent on the initial data by \cite{BS,GLY}. The main difficulty lies in that the DP equation is of hyperbolic type. Precisely speaking, if $u,v\in \mathcal{C}([0,T],B^{1}_{\infty,1})$ are two solutions of \eqref{dp} associated with $u_0$ and $v_0$, in view of \eqref{y3}, we have to tackle with the term $\|v_x\|_{B^1_{\infty,1}}$. To bypass this, we can take $v=\mathbf{S}_{t}(S_Nu_0)$ as the solution to \eqref{dp} with initial data $S_Nu_0$.

Letting $u=\mathbf{S}_{t}(u_0)$ and $v=\mathbf{S}_{t}(S_Nu_0)$, using Lemma \ref{ley1} and \eqref{u2}, we have
\bbal
&\|\mathbf{S}_{t}(S_Nu_0)\|_{B^2_{\infty,1}}\leq C\|S_Nu_0\|_{B^2_{\infty,1}}\leq C2^N\|u_0\|_{B^1_{\infty,1}}
\end{align*}
and
\bbal
\|\mathbf{S}_{t}(S_Nu_0)-\mathbf{S}_{t}(u_0)\|_{L^\infty}&\leq C\|S_Nu_0-u_0\|_{L^\infty}
\\&\leq C\|S_Nu_0-u_0\|_{B^0_{\infty,1}}\\
&\leq C2^{-N}\|S_{N}u_0-u_0\|_{B^1_{\infty,1}},
\end{align*}
which implies
\bbal
\|\mathbf{S}_{t}(S_Nu_0)-\mathbf{S}_{t}(u_0)\|_{B^1_{\infty,1}}&\leq C\|S_Nu_0-u_0\|_{B^1_{\infty,1}}.
\end{align*}
Then, using the triangle inequality, we have for $u_0,\widetilde{u}_0\in B^1_{\infty,1}$,
\bbal
\|\mathbf{S}_{t}(u_0)-\mathbf{S}_{t}(\widetilde{u}_0)\|_{B^1_{\infty,1}}
&\leq \|\mathbf{S}_{t}(S_Nu_0)-\mathbf{S}_{t}(u_0)\|_{B^1_{\infty,1}}+\|\mathbf{S}_{t}(S_N\widetilde{u}_0)-\mathbf{S}_{t}(\widetilde{u}_0)\|_{B^1_{\infty,1}}\\
&~~+\|\mathbf{S}_{t}(S_Nu_0)-\mathbf{S}_{t}(S_N\widetilde{u}_0)\|_{B^1_{\infty,1}}
\\&\leq C\|S_Nu_0-u_0\|_{B^1_{\infty,1}}+C\|S_N\widetilde{u}_0-\widetilde{u}_0\|_{B^1_{\infty,1}}\\
&~~+C\|\mathbf{S}_{t}(S_Nu_0)-\mathbf{S}_{t}(S_N\widetilde{u}_0)\|_{B^1_{\infty,1}}
\\
&:=\mathbf{I}_1+\mathbf{I}_2+\mathbf{I}_3.
\end{align*}
By the interpolation inequality, one has
\bbal
\mathbf{I}_3&\leq C\|\mathbf{S}_{t}(S_Nu_0)-\mathbf{S}_{t}(S_N\widetilde{u}_0)\|^{\fr12}_{B^0_{\infty,\infty}}
\|\mathbf{S}_{t}(S_Nu_0)-\mathbf{S}_{t}(S_N\widetilde{u}_0)\|^{\fr12}_{B^2_{\infty,\infty}}\\
&\leq C\|\mathbf{S}_{t}(S_Nu_0)-\mathbf{S}_{t}(S_N\widetilde{u}_0)\|^{\fr12}_{L^\infty}
\|\mathbf{S}_{t}(S_Nu_0)-\mathbf{S}_{t}(S_N\widetilde{u}_0)\|^{\fr12}_{B^2_{\infty,1}}\\
&\leq C\|S_Nu_0-S_N\widetilde{u}_0\|^{\fr12}_{L^\infty}
\|\mathbf{S}_{t}(S_Nu_0)-\mathbf{S}_{t}(S_N\widetilde{u}_0)\|^{\fr12}_{B^2_{\infty,1}}\\
&\leq C2^{\frac N2}\|u_0-\widetilde{u}_0\|^{\frac12}_{B^1_{\infty,1}},
\end{align*}
which clearly implies
\bbal
\|\mathbf{S}_{t}(u_0)-\mathbf{S}_{t}
(\widetilde{u}_0)\|_{B^1_{\infty,1}}
&\lesssim \|S_Nu_0-u_0\|_{B^1_{\infty,1}}+
\|S_N\widetilde{u}_0-\widetilde{u}_0\|_{B^1_{\infty,1}}+2^{\frac N2}\|u_0-\widetilde{u}_0\|^{\frac12}_{B^1_{\infty,1}}
\\&\lesssim \|S_Nu_0-u_0\|_{B^1_{\infty,1}}+
\|u_0-\widetilde{u}_0\|_{B^1_{\infty,1}}+2^{\frac N2}\|u_0-\widetilde{u}_0\|^{\frac12}_{B^1_{\infty,1}}.
\end{align*}
This completes the proof of continuous dependence.

\subsection{Non-uniform continuous dependence}
 Let $\hat{\phi}\in \mathcal{C}^\infty_0(\mathbb{R})$ be an even, real-valued and non-negative function on $\R$ and satisfy
\begin{numcases}{\hat{\phi}(x)=}
1, &if $|x|\leq \frac{1}{4}$,\nonumber\\
0, &if $|x|\geq \frac{1}{2}$.\nonumber
\end{numcases}
\begin{lemma}\label{ley2}
We define the high frequency function $f_n$ and the low frequency functions $g_n$ as follows
\bbal
&f_n=2^{-n}\phi(x)\sin \bi(\frac{17}{12}2^nx\bi),\\
&g_n=\frac{12}{17}2^{-n}\phi(x),\quad n\gg1.
\end{align*}
Then for any $\sigma\in\R$, we have
\bal
&\|f_n\|_{L^\infty}\leq C2^{-n}\phi(0)\quad\text{and}\quad\|g_n\|_{L^\infty}\leq C2^{-n}\phi(0),\label{y0}\\
&\|f_n\|_{B^{\sigma}_{\infty,1}}\leq C2^{(\sigma-1)n}\phi(0)\quad\text{and}\quad\|g_n\|_{B^\sigma_{\infty,1}}\leq C2^{-(n+\sigma)}\phi(0),\label{yz}\\
&\liminf_{n\rightarrow \infty}\|g_n\pa_xf_n\|_{B^{1}_{\infty,\infty}}\geq M_1,\label{yz3}
\end{align}
for some positive constants $C, M_1$.
\end{lemma}
{\bf Proof.}\quad We refer to see Lemma 3.1 in \cite{Li2} for the proof with minor modifications.
\begin{proposition}\label{pro1}
Assume that $\|u_0\|_{B^{1}_{\infty,1}}\lesssim 1$. Under the assumptions of Theorem \ref{th1}, we have
\bal\label{et0}
\|\mathbf{S}_{t}(u_0)-u_0-t\mathbf{v}_0(u_0)\|_{B^{1}_{\infty,1}}\leq Ct^{2}\mathbf{E}(u_0),
\end{align}
where we denote $\mathbf{v}_0(u_0):=\mathbf{P}(u_0)-u_0\pa_x u_0$ and
\bbal
\mathbf{E}(u_0)&:=1+\|u_0\|_{L^\infty}\big(\|u_0\|_{B^{2}_{\infty,1}}+
\|u_0\|_{L^\infty}
\|u_0\|_{B^{3}_{\infty,1}}\big).
\end{align*}
\end{proposition}
{\bf Proof.}\quad For simplicity, we denote $u(t)=\mathbf{S}_t(u_0)$.
By the Mean Value Theorem, we obtain
\bal\label{et1}
\|u(t)-u_0\|_{L^\infty}&\leq\int^t_0\|\pa_\tau u\|_{L^\infty} \dd\tau
\nonumber\\&\leq \int^t_0\|u\pa_xu\|_{L^\infty}\dd \tau+\int^t_0\|\mathbf{P}(u)\|_{L^\infty}\dd \tau
\nonumber\\&\leq C\int^t_0\|u\|_{L^\infty}\|u\|_{C^{0,1}}\dd \tau\nonumber\\
&\leq Ct\|u_0\|_{L^\infty},
\end{align}
where we have used the estimate
$$\|\mathbf{P}\big(u\big)\|_{L^\infty}\leq C\|u\|^2_{L^\infty}.$$
Using Lemma \ref{le2}, \eqref{u2}  and  \eqref{et1} yield
\bal\label{et2}
\|u(t)-u_0\|_{B^{1}_{\infty,1}}
&\leq \int^t_0\|\pa_\tau u\|_{B^{1}_{\infty,1}} \dd\tau
\nonumber\\&\leq \int^t_0\|\mathbf{P}(u)\|_{B^{1}_{\infty,1}} \dd\tau+ \int^t_0\|u \pa_xu\|_{B^{1}_{\infty,1}} \dd\tau
\nonumber\\&\leq Ct\big(\|u\|^{2}_{B^{1}_{\infty,1}}+\|u\|_{L^\infty}\|u_x\|_{B^{1}_{\infty,1}}\big)
\nonumber\\&\leq Ct\big(1+\|u_0\|_{L^\infty}
\|u_0\|_{B^{2}_{\infty,1}}\big).
\end{align}
Similarly, we have
\bal\label{et3}
\|u(t)-u_0\|_{B^{2}_{\infty,1}}
&\leq \int^t_0\|\pa_\tau u\|_{B^{2}_{\infty,1}} \dd\tau
\nonumber\\&\leq \int^t_0\|\mathbf{P}(u)\|_{B^{2}_{\infty,1}} \dd\tau+ \int^t_0\|u \pa_xu\|_{B^{2}_{\infty,1}} \dd\tau
\nonumber\\&\leq Ct\big(\|u\|_{B^{1}_{\infty,1}}\|u\|_{B^{2}_{\infty,1}}
+\|u\|_{L^\infty}\|u\|_{B^{3}_{\infty,1}}\big)
\nonumber\\&\leq Ct\big(\|u_0\|_{B^{2}_{\infty,1}}
+\|u_0\|_{L^\infty}\|u_0\|_{B^{3}_{\infty,1}}\big).
\end{align}
Using the Mean Value Theorem and Lemma \ref{le2} once again, we obtain that
\bal\label{et4}
\|u(t)-u_0-t\mathbf{v}_0(u_0)\|_{B^{1}_{\infty,1}}
&\leq \int^t_0\|\pa_\tau u-\mathbf{v}_0(u_0)\|_{B^{1}_{\infty,1}} \dd\tau
\nonumber\\&\leq \int^t_0\|\mathbf{P}(u)-\mathbf{P}(u_0)\|_{B^{1}_{\infty,1}} \dd\tau+\int^t_0\|u\pa_xu-u_0\pa_xu_0\|_{B^{1}_{\infty,1}} \dd\tau
\nonumber\\&\lesssim \int^t_0\|u(\tau)-u_0\|_{B^{1}_{\infty,1}} \dd\tau+\int^t_0\|u(\tau)-u_0\|_{L^\infty} \|u(\tau)\|_{B^{2}_{\infty,1}} \dd\tau
\nonumber\\&\quad \ + \int^t_0\|u(\tau)-u_0\|_{B^{2}_{\infty,1}}  \|u_0\|_{L^\infty}\dd \tau\nonumber\\
&\lesssim \int^t_0\|u(\tau)-u_0\|_{B^{1}_{\infty,1}} \dd\tau+\|u_0\|_{B^{2}_{\infty,1}}\int^t_0\|u(\tau)-u_0\|_{L^\infty}  \dd\tau
\nonumber\\&\quad \ + \|u_0\|_{L^\infty}\int^t_0\|u(\tau)-u_0\|_{B^{2}_{\infty,1}}  \dd \tau.
\end{align}
Plugging \eqref{et1}--\eqref{et3} into \eqref{et4} yields the desired result \eqref{et0}. Thus, we complete the proof of Proposition \ref{pro1}.

Now we prove the non-uniform continuous dependence.

We set $u^n_0=f_n+g_n$ and compare the solution $\mathbf{S}_{t}(u^n_0)$ and $\mathbf{S}_{t}(f_n)$. We obviously have
\bbal
\|u^n_0-f_n\|_{B^{1}_{\infty,1}}=\|g_n\|_{B^{1}_{\infty,1}}\leq C2^{-n},
\end{align*}
which means that
\bbal
\lim_{n\to\infty}\|u^n_0-f_n\|_{B^{1}_{\infty,1}}=0.
\end{align*}
From Lemma \ref{le2}, one has
\bbal
&\|u^n_0,f_n\|_{B^{\sigma}_{\infty,1}}\leq C2^{(\sigma-1)n}\quad \text{for}\quad \sigma\geq1,\\
&\|u^n_0,f_n\|_{L^\infty}\leq C2^{-n},
\end{align*}
which implies
\bbal
\mathbf{E}(u^n_0)+\mathbf{E}(f_n)\leq C.
\end{align*}
Using the facts
\bbal
&\big\|u^n_{0}\pa_xg_n\big\|_{B^{1}_{\infty,1}}\leq C\big\|u^n_0\big\|_{B^{1}_{\infty,1}}\big\|g_n\big\|_{B^{2}_{\infty,1}}\leq C2^{-n},\\
&\big\|\mathbf{P}(u^n_0)-\mathbf{P}(f_n)\big\|_{B^{1}_{\infty,1}}\leq C\big\|g_n\big\|_{B^{1}_{\infty,1}}\big\|u^n_0+f_n\big\|_{B^{1}_{\infty,1}}\leq C2^{-n}.
\end{align*}
we deduce that (for more details, see \cite{Li2})
\bal\label{yyh}
\big\|\mathbf{S}_{t}(u^n_0)-\mathbf{S}_{t}(f_n)\big\|_{B^{1}_{\infty,1}}\geq&~ t\big\|g_n\pa_xf_n\big\|_{B^{1}_{\infty,1}}-t\big\|u^n_{0}\pa_xg_n,\;\mathbf{P}(u^n_0)-\mathbf{P}(f_n)\big\|_{B^{1}_{\infty,1}}-Ct^{2}-C2^{-n}\nonumber\\
\geq&~ t\big\|g_n\pa_xf_n\big\|_{B^{1}_{\infty,1}}-Ct2^{-n}-Ct^{2}-C2^{-n},
\end{align}
Notice that \eqref{yz3}
\begin{eqnarray*}
      \liminf_{n\rightarrow \infty} \big\|g_n\pa_xf_n\big\|_{B^{1}_{\infty,1}}\gtrsim M_1,
        \end{eqnarray*}
then we deduce from \eqref{yyh} that
\bbal
\liminf_{n\rightarrow \infty}\big\|\mathbf{S}_t(f_n+g_n)-\mathbf{S}_t(f_n)\big\|_{B^{1}_{\infty,1}}\gtrsim t\quad\text{for} \ t \ \text{small enough}.
\end{align*}
This completes the proof of Theorem \ref{th1}.

\section*{Acknowledgments} J. Li is supported by the National Natural Science Foundation of China (Grant No.11801090). Y. Yu is supported by the Natural Science Foundation of Anhui Province (No.1908085QA05). W. Zhu is partially supported by the National Natural Science Foundation of China (Grant No.11901092) and Natural Science Foundation of Guangdong Province (No.2017A030310634).


\begin{thebibliography}{99}
\linespread{0}\addtolength{\itemsep}{-1.0ex}

\bibitem{B.C.D} H. Bahouri, J. Y. Chemin, R. Danchin, Fourier Analysis and Nonlinear Partial Differential Equations, Grundlehren der Mathematischen Wissenschaften, Springer, Heidelberg, 2011.

\bibitem{BS} J. L. Bona and R. Smith, The initial-value problem for the Korteweg-de Vries equation, Philos. Trans. R. Soc. Lond. Ser. A, 278 (1975), 555-601.

\bibitem{Camassa} R. Camassa, D. Holm, An integrable shallow water equation with peaked solitons, Phys. Rev. Lett., 71 (1993), 1661-1664.
\bibitem{Dullin} H.R. Dullin, G.A. Gottwald, D.D. Holm, On asymptotically equivalent shallow water wave equations, Physica D, 190 (2004), 1-14.
\bibitem{DP} A. Degasperis, D. Holm, A. Hone, A new integral equation with peakon solutions. Theoret. Math.
Phys. 133 (2002), 1463-1474.

\bibitem{DP1} A. Degasperis, M. Procesi, Asymptotic integrability. In: Symmetry and Perturbation Theory. (Rome, 1998). Rivers Edge, NJ: World Scientific
Publishing Company, pp. 23-37 (1999).

\bibitem{Escher} J. Escher, Y. Liu, Z. Yin, Shock waves and blow-up phenomena for the periodic Degasperis-Procesi equation, Indiana
Univ. Math. J. 56 (2007), 87-177.

\bibitem{Escherjfa} J. Escher, Y. Liu, Z. Yin, Global weak solutions and blow-up structure for the Degasperis-Procesi equation, J. Funct.
Anal. 241 (2006), 457-485.

\bibitem{Gui} G. Gui, Y. Liu, On the Cauchy problem for the Degasperis-Procesi equation, Quart. Appl. Math. 69 (2011), 445-464.

\bibitem{GLY} Z. Guo, J. Li, Z. Yin, Local well-posedness of the incompressible Euler equations in $B_{\infty,1}^1$ and the inviscid limit of the Navier-Stokes equations, J. Funct. Anal., 276 (2019), 2821--2830.

\bibitem{Guo-Yin} Z. Guo, X. Liu, L. Molinet, Z. Yin, Ill-posedness of the Camassa-Holm and related equations in the critical space, J. Differential Equations, 266 (2019), 1698-1707.



\bibitem{Himonas2014} A. Himonas, C. Holliman, K. Grayshan, Norm Inflation and Ill-Posedness for the Degasperis-Procesi Equation, Comm. Partial Differ. Equ. 39 (2014), 2198-2215.

\bibitem{Himonas2016} A. Himonas, K. Grayshan, C. Holliman, Ill-Posedness for the b-Family of Equations, J. Nonlinear Sci. 26 (2016), 1175-1190.


\bibitem{Himonas} A.A. Himonas, C. Holliman, The Cauchy problem for the Novikov equation, Nonlinearity 25 (2012) 449-479.

\bibitem{Lundmark} H. Lundmark, J. Szmigielski, Multi-peakon solutions of the Degasperis-Procesi equation, Inverse Problems. 19 (2003),
1241-1245.
\bibitem{Lundmark2007} H. Lundmark, Formation and dynamics of shock waves in the Degasperis-Procesi equation, J. Nonlinear Sci. 17 (2007),
169-198.
\bibitem{Lenells} J. Lenells, Traveling wave solutions of the Degasperis-Procesi equation, J. Math. Anal. Appl. 306 (2005), 72-82.
\bibitem{Li0} J. Li, Z. Yin, Remarks on the well-posedness of Camassa-Holm type equations in Besov spaces, J. Differential Equations, 261 (2016), 6125-6143.

\bibitem{Li1} J. Li, Y. Yu, W. Zhu, Non-uniform dependence on initial data for the Camassa-Holm equation in Besov spaces, J. Differential Equations, 269 (2020), 8686--8700.

\bibitem{Li2} J. Li, X. Wu, Y. Yu, W. Zhu, Non-uniform dependence on initial data for the Camassa-Holm equation in the critical Besov space. J. Math. Fluid Mech., 23:36 (2021), 11 pp.
\bibitem{Liu2006} Y. Liu, Z. Yin, Global existence and blow-up phenomena for the Degasperis-Procesi equation, Commun. Math. Phys. 267 (2006), 801-820.
\bibitem{Liuyin} Y. Liu, Z. Yin, On the blow-up phenomena for the Degasperis-Procesi equation, Int. Math. Res. Not. IMRN 23 (2007),
117. 22 pp.
\bibitem{Lin} Z. Lin, Y. Liu, Stability of peakons for the Degasperis-Procesi equation, Commun. Pure Appl. Math. 62(1) (2009), 125-146.
\bibitem{Vakhnenko} V.O. Vakhnenko, E.J. Parkes, Periodic and solitary-wave solutions of the Degasperis-Procesi equation, Chaos Solitons
Fractals. 20 (2004), 1059-1073.
\bibitem{Yin2003} Z. Yin,  Global existence for a new periodic integrable equation, J. Math. Anal. Appl. 283(2003), 129-139.
\bibitem{Yin2004} Z. Yin, Global solutions to a new integrable equation with peakons, Indiana Univ. Math. J. 53 (2004), 1189-1210.
\bibitem{Yin2} Z. Yin, On the Cauchy problem for an integrable equation with peakon
solutions, Illinois J. Math. 47, 649-666 (2003).
\bibitem{YinJFA} Z. Yin, Global weak solutions for a new periodic integrable equation with peakon solutions, J. Funct. Anal. 212 (2004),
182-194.




\end{thebibliography}
\end{document}